\documentclass[12pt,reqno]{amsart}

\usepackage[dvips]{graphicx} 

\usepackage{%amssymb, latexsym, amsmath, amscd, array, 
hyperref
}

\usepackage{breakurl}   %this allows for url breaks in pdf

\usepackage{comment} %for longer comments

\usepackage{apalike} %to change citations to parentheses instead of
\theoremstyle{definition}

\numberwithin{equation}{section}

\newcommand\N {{\mathbb N}} 

\newcommand\R {{\mathbb R}}

\renewenvironment{quote}
               {\list{}{}%
                \item\relax}
               {\endlist}

\title{Formalism 25}

\author[M. Katz]{Mikhail G. Katz}\address{Department of Mathematics,
  Bar Ilan University, Ramat Gan 5290002 Israel
  \url{http://orcid.org/0000-0002-3489-0158}}
\email{katzmik@math.biu.ac.il}

\author[K. Kuhlemann]{Karl Kuhlemann}\address{Gottfried Wilhelm
  Leibniz University Hannover, D-30167 Hannover, Germany
  \url{http://orcid.org/0000-0002-7713-4782}}
\email{kus.kuhlemann@t-online.de}

\author[S. Sanders]{Sam Sanders} \address{Department of Philosophy 2,
  RUB Bochum, Bochum, Germany \url{http://sasander.wix.com/academic}
  \url{https://orcid.org/0000-0001-8256-0009}}
\email{sasander@me.com}

\author[D. Sherry]{David Sherry} \address{Department of Philosophy,
  Northern Arizona University, Flagstaff, AZ 86011, US
  \url{http://orcid.org/0000-0001-9699-7762}}
\email{David.Sherry@nau.edu}

\subjclass[2020]{Primary 00A30 %phil of math
Secondary 01A60, %history and biography, 20th century
03A05
}

\begin{document}

\begin{abstract}
Abraham Robinson's philosophical stance has been the subject of
several recent studies.  Erhardt following Gaifman claims that
Robinson was a finitist, and that there is a tension between his
philosophical position and his actual mathematical output.  We present
evidence in Robinson's writing that he is more accurately described as
adhering to the philosophical approach of Formalism.  Furthermore, we
show that Robinson explicitly argued \emph{against} certain finitist
positions in his philosophical writings.  There is no tension between
Robinson's mathematical work and his philosophy because mathematics
and metamathematics are distinct fields: Robinson advocates finitism
for metamathematics but no such restriction for mathematics.  We show
that Erhardt's analysis is marred by historical errors, by routine
conflation of the generic and the technical meaning of several key
terms, and by a philosophical \emph{parti pris}.  Robinson's Formalism
remains a viable alternative to mathematical Platonism.
\end{abstract}

\today

%\doublespacing
\thispagestyle{empty}

%\huge

\keywords{Finitism; Formalism; foundations; Platonism; Realism;
  Robinson}

\maketitle
%\tableofcontents

\section{Formalism and finitism}

Abraham Robinson's intellectual profile continues to inspire passions
fifty years after his passing.  A number of publications over the past
few decades have analyzed his philosophical position, including
\cite[p.\;603]{Be88}, \cite{Ga12}, \cite{14c}, \cite{17i},
\cite{We24}, and most recently \cite{Er25}.  A~majority of
commentators take it for granted that Robinson's mature position is
best described as mathematical Formalism.  However, Erhardt reads
Robinson as a finitist (see Section~\ref{s7} for a discussion of
distinct meanings of the term) who rejects infinite totalities.

\subsection{Meaning and reference}
\label{s1}

Significantly, there is a key observation in Robinson's 1975 article
that is missing from Erhardt's text (though he does cite the article).
The observation helps understand Robinson's take on \emph{meaning}.
In 1975, Robinson clarified his position by stating that
\begin{quote}
mathematical theories that, allegedly, deal with infinite totalities
have no detailed meaning, i.e.\;\emph{reference}.  \cite[p.\;42;
  emphasis added]{Ro75}; reprinted in \cite[p.\;557]{Ro79}
\end{quote}
Namely, the term \emph{infinite totality} has no reference (or
referent) in either the physical world or any Platonic realm of
mathematical abstracta.  Robinson's main goal here was to distance
himself from mathematical Platonism.  He did \emph{not} believe that
expressions such as \emph{infinite totality} lacked meaning in the
sense of being `pointless' or `devoid of significance', as we will see
below; he only intended that such expressions lacked a reference, as
we explained.  In his earlier philosophical text \emph{Formalism\;64}
he used the term \emph{direct interpretation} in place of
\emph{reference}:
\begin{quote}
I {\ldots}~regard a theory which refers to an infinite totality as
\emph{meaningless} in the sense that its terms and sentences cannot
possess the \emph{direct interpretation} in an actual structure that we
should expect them to have by analogy with concrete (e.g., empirical)
situations.%
\footnote{\cite[p.\;231]{Ro65}; emphasis on `meaningless' in
the original; emphasis on `direct interpretation' added.}
\end{quote}
Having outlined Robinson's position, let us examine Erhardt's
presentation thereof.  Erhardt misrepresents Robinson's position by
quoting Robinson out of context.  Thus, Erhardt claims the following:
\begin{quote}
In \emph{Formalism 64}, Robinson’s most fully developed philosophical
paper, he writes (Robinson 1964/1979b, pp.\;230–231, italics in the
original), ``[\ldots] the notion of a particular \emph{class} of five
elements, e.g. of five particular chairs, presents itself to my mind
as clearly as the notion of a single individual [\ldots].  By
contrast, I feel quite \emph{unable to grasp} the idea of an actual
infinite totality.  To me there appears to exist an unbridgeable gulf
between sets or structures of one, or two, or five elements, on one
hand, and infinite structures on the other hand~[\ldots]''%
\footnote{Robinson as quoted in \cite[p.\;431]{Er25}; emphasis on
``unable to grasp'' added.}
\end{quote}
We stress the profusion of ellipses `[\ldots]' in Erhardt's
quotation from Robinson.  Superficially, the abridged passage may
suggest a finitist position.  It is instructive to explore the issue
whether Erhardt's deletions alter the meaning of the passage as
intended by Robinson.  On the face of it, the passage as quoted by
Erhardt sounds rather odd.  A modern mathematician who claims to be
``unable to grasp'' the concept of an actual infinite totality
certainly sounds peculiar.  We will analyze Robinson's passage in its
context in Section~\ref{s2}.

\subsection{Gaifman's reading of Robinson}
\label{s7}

Erhardt thanks Gaifman for encouraging him to read Robinson's text
\emph{Formalism 64} {\cite[p.\;446]{Er25}}.
Gaifman may have been the original source of a characterisation of
Robinson as a finitist.  Indeed, Gaifman claims in 2012:
\begin{quote}
Abraham Robinson, who was a \emph{finitist}, or something very near to
it, realized the seriousness of the limitations that his position
implied with regard to syntactic concepts that required quantification
over infinite domains. {\cite[p.\;488; emphasis added]{Ga12}}
\end{quote}
Gaifman's own position is analyzed in Section~\ref{s12}.  Was Robinson
a finitist as claimed by Gaifman, and does Robinson's position entail
serious `limitations' as Gaifman claims?  To address the issue, it is
crucial to distinguish between finitism in a narrow sense and finitism
in a broad sense.  In its broad sense, finitism denotes opposition to
the use of infinitary concepts at the metamathematical level (as in
Hilbert's program).  In its narrow sense, finitism is characterized by
an opposition to the use of infinite totalities at both the
metamathematical and the mathematical level.  Thus, many intuitionists
and constructivists were opposed to the use of \emph{certain} infinite
totalities in mathematical practice.

Erhardt's claim of tension in Robinson's work stems from a conflation
of these two meanings of finitism.  The fact that Robinson's position
is more accurately described as Formalism than finitism is evident
from his recommendation concerning the business of mathematics:
\begin{quote}
[W]e should continue the business of Mathematics ``as usual,'' i.e.,
we should \emph{act as if infinite totalities really existed}.
\\ {\cite[p.\;230; emphasis added]{Ro65}}
\end{quote}
Finitists (in the narrow sense) \emph{never} ``act as if infinite
totalities really existed.''  The idea that Robinson was not a
finitist in the narrow sense is not merely a matter of our opinion
against Erhardt's.  Indeed, Robinson explicitly argued \emph{against}
some finitists' rejection of the use of infinitary terms in
mathematics:
\begin{quote}
Those who adopt this attitude [including the Intuitionists] think that
a concept, or a sentence, or an entire theory, is acceptable only if
it can be \emph{understood} properly and that a concept, or sentence,
or a theory, is understood properly only if all terms which occur in
it can be interpreted directly, as explained. By contrast, the
formalist holds that direct interpretability is not a necessary
condition for the acceptability of a mathematical theory.
{\cite[p.\;234; emphasis in the original]{Ro65}}
\end{quote}
One such intuitionist is Dummett, who claimed the following:
\begin{quote}
Constructivist philosophies of mathematics insist that the
\emph{meanings of all terms}, including logical constants, appearing
in mathematical statements must be given in relation to constructions
which we are capable of effecting, and of our capacity to recognise
such constructions as providing proofs of those statements; {\ldots}
{\cite[p.\;301; emphasis added]{Du75}}
\end{quote}
Robinson clearly disagreed with Dummett's claim that \emph{all terms}
must be assigned such a direct meaning.  Robinson concluded:
\begin{quote}
To sum up, the direct interpretability of the terms of a mathematical
theory is not a necessary condition for its acceptability; a theory
which includes \emph{infinitary terms} is not thereby less acceptable
or less rational than a theory which avoids them.  {\cite[p.\;235;
    emphasis added]{Ro65}}
\end{quote}
Here Robinson endorsed the Formalist position and moreover contrasted
it with finitism in the narrow sense.
%Erhardt's text contains further claims but they all rest on his
%dubious assumption that Robinson was a finitist.

\subsection{Standard model}
\label{s13}

The source of Platonists' discomfort with Formalism (and their
proclivity to paint Formalists as finitists) is identified in section
10 of Robinson's text \emph{Formalism 64}:
\begin{quote}
In particular, I will mention here the assumption that there exists a
\emph{standard} or \emph{intended model} of Arithmetic or
(alternatively, but relatedly) of Set Theory.  Clearly, to the
formalist, the entire notion of standardness must be meaningless, in
accordance with our first basic principle.%
\footnote{\cite[p.\;242]{Ro65}; emphasis in the original.  For a
discussion of this use of the term \emph{standardness} see
Section~\ref{s8}.}
\end{quote}
As discussed in Section~\ref{s1}, Robinson denied the existence of a
referent for such a standard model (a.k.a.~intended interpretation) of
$\mathbb N$ or~$\mathbb R$, in either the physical or any Platonic
realm.  To Platonists, this may appear as a narrow finitist stance,
but they may well ponder why Robinson did not name his article
``Finitism 64''.

Erhardt goes on to comment as follows:
\begin{quote}
There is, however, another alternative: more strongly committing to
finitism. This is the path taken by various \emph{constructivists},
who rather than play the uninterpretable game of symbols choose, in a
philosophically principled manner, to adopt a weaker form of
mathematics.  \cite[p.\,444; emphasis added]{Er25}
\end{quote}
Here Erhardt may be right about constructivists (rather than Robinson)
describing classical mathematics as an `uninterpretable game of
symbols',%
\footnote{\label{f11b}Thus, Bishop writes: ``The successful
formalization of mathematics helped keep mathematics on a wrong
course.
%The fact that space has been arithmetized loses much of its
%significance if space, number, and everything else are fitted into a
%matrix of idealism where even the positive integers have an ambiguous
%computational existence.
{\ldots}~Mathematics becomes the game of sets, which is a fine game as
far as it goes, with rules that are admirably precise''
\cite[p.\;4]{Bi67}.  In the same vein: ``If every mathematician
occasionally, perhaps only for an instant, feels an urge to move
closer to reality, it is not because he believes mathematics is
lacking in meaning. He does not believe that mathematics consists in
drawing brilliant conclusions from arbitrary axioms, of juggling
concepts devoid of pragmatic content, of playing a meaningless game''
\cite[p.\;viii]{Bi67}.}
though Erhardt does not admit in his article that Robinson himself
does \emph{not} describe it in these terms.  Erhardt's first sentence
is based on the unfounded assumption that Robinson was a finitist, as
we have already discussed; Erhardt compounds his error - of
attributing to Robinson the position of narrow finitism%
\footnote{As is evident from Erhardt's comment on Robinson's would-be
reaction to Wiles' proof of Fermat's Last Theorem; see
Section~\ref{Fermat}.}
- by taking it upon himself to offer advice as to how to
practice the latter.

\subsection{Robinson's passage in context}
\label{s2}

In Section~\ref{s1}, we saw that Erhardt claimed that Robinson was a
finitist based on a passage from Robinson's text \emph{Formalism 64}.
Examining Robinson's passage in context reveals a rather different
picture.  Robinson begins by mentioning the traditional positions of
mathematical philosophy (Formalism, Intuitionism, Logicism) on
page~228 of his text \emph{Formalism 64}.  On page~230, he mentions a
philosophical school of nominalism.%
\footnote{Today the term \emph{nominalism} may refer to any
anti-Platonist philosophy of mathematics.  Robinson uses the term in a
narrower sense.}
Robinson does not treat nominalism in much detail, on the grounds that
to him it represents ``little depth from the mathematical point of
view'' {\cite[p.\;231]{Ro65}}.
Robinson claims that to a nominalist, there is not much difference
between a set of five objects and an infinite set.  Both are
`illusory' to a nominalist, according to Robinson:
\begin{quote}
To a nominalist, the existence of a set of five elements is no less
\emph{illusory} than the existence of the totality of all natural
numbers.  At the other end of the scale are the so-called platonic
realists or platonists who believe in the \emph{ideal} existence of
mathematical entities in general, including the existence of
transfinite sets of arbitrarily large cardinal numbers to the extent
to which they can be introduced at all by means of suitable axioms.
{\cite[p.\;230; emphasis added]{Ro65}}
\end{quote}
By the time Robinson gets to the ``five particular chairs'' (as quoted
by Erhardt) at the bottom of page 230, it is clear that his goal is to
stress the difference between his position and that of nominalists.
After noting the difference between ``five particular chairs'' and an
infinite set, Robinson emphasizes on page 231 that describing infinite
totalities as `meaningless' does \emph{not} mean that ``such a theory
is therefore pointless or devoid of significance''
{\cite[p.\;231]{Ro65}.}

As noted in Section~\ref{s1}, Robinson clarified his adjective
`meaningless' in his 1975 text where he stated that what he has in
mind is an absence of a \emph{reference}; namely, the term
\emph{infinite totality} does not \emph{refer} to anything in either
the real or a Platonic realm.  Already in \emph{Formalism\;64} he
emphasized that he rejected the idea that infinite totalities exist
``either really or ideally''  {\cite[p.\;230]{Ro65}.}
Here `really' refers to the physical realm, and `ideally' to a
Platonic realm, as he states explicitly in the passage quoted above.
It emerges that Erhardt's presentation of Robinson's position, as when
Erhardt claims that it ``places Robinson at odds with typical
mathematical practice'' \cite[p.\;432]{Er25},
conflates the \emph{generic} meaning of terms such as \emph{to~grasp}
and \emph{meaningless}, with the precise \emph{technical} meaning
attributed to such terms by Robinson.

\subsection{Fermat's Last Theorem: is the proof legitimate?}
\label{Fermat}

The conflation of the generic and the technical meaning of the term
\emph{meaningless} has a further effect of leading Erhardt to a
preposterous misrepresentation of Robinson's position, as when Erhardt
claims:
\begin{quote}
Though Fermat’s Last Theorem is a general result proven by
illegitimate methods and that purports to say something about all
natural numbers—constituting an illicit reference to actual
infinity—one can derive from it a potentially infinite number of
legitimate, material claims. Yet the fact that the proof is
\emph{illegitimate on Robinson’s view} is crucial to our appraisal of
his view.  \cite[p.\;441; emphasis added]{Er25}
\end{quote}
Erhardt's \emph{illegitimacy} claim stems from his mistaken
identification of Robinson as a narrow finitist.  But Wiles' proof of
Fermat's Last Theorem would \emph{not} be `illegitimate on Robinson's
view' as Erhardt claims.  While Wiles' proof may seem illegitimate to
those finitists who view proofs involving infinite totalities as
\emph{meaningless} in the generic sense of the term, it is certainly
legitimate on Robinson's view, collapsing what Erhardt describes as
``our appraisal of his view'' due to Erhardt's conflation of distinct
meanings of the term \emph{meaningless}.

While discussing Robinson's commitment to potential infinity, Erhardt
claims that Robinson is a `devout realist' about each natural number:
\begin{quote}
In opposition to nominalists (and \emph{a fortiori}, fictionalists),
Robinson is a \emph{devout realist} with regard to each individual
natural number.  This is evidenced by his commitment to potential
infinity, {\ldots}%
\footnote{\cite[p.\;438, note 22]{Er25}; emphasis on ``devout
realist'' added.}
\end{quote}
Robinson himself would have likely rejected such a description of his
position.  Being able to grasp a concept (such as a potential infinity
of metalanguage integers) does not commit one to a reality of a
Platonic object.
%One cannot be a realist about something that is only potentially
%infinite.
It emerges that Erhardt is trying to foist an untenable position on
Robinson.

With regard to the issue of nominalism and fictionalism mentioned by
Erhardt above, we note a pecularity of Erhardt's approach.  Erhardt is
keen to cast Robinson as an instrumentalist, rather than one who
regards infinite classes as useful fictions.  Erhardt dismisses the
suggestion that Robinson was a fictionalist, in spite of acknowledging
that Leibniz, an avowed fictionalist, ``had a profound effect'' on
Robinson \cite[p.\;437 and note 22]{Er25}.
The rationale for such a dismissal is most peculiar: Erhardt considers
that it would be anachronistic for Robinson, in 1964, to see himself
as a fictionalist since Field's `introduction of fictionalism'
occurred only in 1980 (see \cite{Fi80}).  Erhardt is right to note
that Field's fictionalism, a form of nominalism, conflicts with
Robinson's repeated rejection of nominalism.  But such a conflict
arises only if fictionalism is assumed to be necessarily a nominalist
proposal.  Leibniz's fictionalism offers an alternative to nominalism.
Leibniz treated infinitesimals, along with negative and imaginary
numbers, as well-founded fictions; see \cite{14c}.  Although he
sometimes suggested that infinitesimals are eliminable in the manner
of Archimedean exhaustion arguments, he made no such suggestion for
negatives and imaginaries.  In all three cases, it is the contribution
to systematicity that establishes the mantle `well-founded fiction',
rather than any sort of nominalist reduction.

\subsection{Symptomatic keyword list and imaginary tensions}

Already Erhardt's keyword list is symptomatic of a problem with his
text: the keyword \emph{finitism} appears, but \emph{formalism} does
not.  Erhardt misinterprets Robinson's comment about being `unable to
grasp' actual infinite collections, by viewing it as a finitist
stance.  As analyzed in Section~\ref{s2}, Robinson merely sought to
distance himself from nominalism and Platonism.  Erhardt's Abstract
claims to detect a `tension' between Robinson's philosophy and his
mathematical practice:
\begin{quote}
The foundational position he inherited from David Hilbert undermines
not only the use of nonstandard analysis, but also Robinson's
considerable corpus of pre-logic contributions%
\footnote{\label{f11}Erhardt's wording here constitutes a historical
error.  He makes it appear as though Robinson's work in aeronautics,
such as his 1956 book \emph{Wing Theory}~\cite{Ro56}, \emph{preceded}
his work in logic.  This is incorrect.  For example, Robinson's
dissertation \emph{The metamathematics of algebraic systems}
\cite{Ro49}, advised by Paul Dienes, dates from 1949; see
\url{https://www.genealogy.math.ndsu.nodak.edu/id.php?id=15886}.  See
also Dauben \cite[p.\,157]{Da95}.}
to the field in such diverse areas as differential equations and
aeronautics.  This \emph{tension} emerges from Robinson's disbelief in
the existence of infinite totalities {\ldots} {\cite[Astract; emphasis
    added]{Er25}}
\end{quote}
Erhardt is referring to the book \emph{Non-standard Analysis}
\cite{Ro66}.  However, Erhardt's argument for his claim is too strong
because it would apply to all Formalists.  Erhardt's Abstract accuses
Robinson of ``giving up on a commitment to reconciling'' his
philosophy and his mathematical practice.  Possibly a Platonist would
tend to view all Formalists in such a fashion.  This may betray
Erhardt's own philosophical stance.
%(see further in Section~\ref{s5}).

In sum, Erhardt's depiction of Robinson's philosophical position falls
prey to a straw man fallacy.  Erhardt's claim in his abstract that
Robinson's philosophical position is somehow at odds with his work in
applied mathematics, is unfounded; see further in Section~\ref{s8}.

\section{Twin prime conjecture and types of infinite totalities}

In his introduction, Erhardt claims that
\begin{quote}
[E]ven basic conjectures of number theory, such as the twin prime
conjecture, presuppose the existence of infinite totalities.
{\cite[p.\;431]{Er25}}
\end{quote}
However, such a claim involves a conflation of distinct meanings
associated with the term \emph{infinite totality}; namely, of
different levels of language or theory.  The twin prime conjecture
(TPC) can be formulated in Peano Arithmetic (PA).  This can be done,
for example, by the formula on page 442, line 4 in \cite{Er25}, which
we will express as follows:
\begin{equation}
  \label{e41}
\forall x \, \exists y \, \phi(x,y),
\end{equation}
where the formula~$\phi$ says that~$y>x$ and both~$y$ and~$y+2$ are
prime.  Note that PA is a theory that does not know any infinite sets
and that is even biinterpretable with the theory obtained from ZF by
replacing the axiom of infinity by its negation, and
adding~$\in$-induction;
%see Caicedo's answer at https://math.stackexchange.com/a/107669/72694
%
{see \cite{Ac37}; \cite{Ka07}.}
%
%The models of PA must of course be infinite.
The claim that the TPC presupposes the existence of infinite
totalities is incorrect with respect to the object theory PA in which
the TPC is formulated (because PA does not know any infinite
totalities).

On the other hand, such a claim is true with respect to a metatheory
in which PA is interpreted in an appropriate structure and thus given
a semantics.
%The domain of such a structure is an infinite totality.
%However, this is “meaningless” for Robinson, since it has no referent.
On this reading, the term \emph{infinite totality} refers to
quantification over an infinite domain.  Since formulation~\eqref{e41}
involves such quantification, it can be said to involve \emph{infinite
totalities} in such a different sense.  The special status of
$\Sigma^0_1$ and $\Pi^0_1$-sentences is discussed in Section~\ref{s12}.

%Erhardt mentions a further principle to the effect that ``we are
%ontologically committed to the objects over which we quantify''%
%
%\footnote{Erhardt \cite[note 2]{Er25}.}
%
%(and attributes such a principle to Quine \cite{Qu48}).  However,
%Erhardt's principle turns out to be too crude a tool to determine the
%epistemological status of mathematical assertions.  Indeed, the
%dependence of TPC on quantification over an infinite domain turns out
%to be an illusion, for the following reverse-mathematical reasons.

%Erhardt's claim that sentences ``presuppose the existence of infinite
%totalities'' does not apply to~$\Pi^0_2$-sentences in full generality;
%namely, it does not apply to sentences that are provable in the theory
%WKL$_0$.  This is because WKL$_0$ is a~$\Pi^0_2$-conservative
%extension of Primitive Recursive Arithmetic, which is generally
%considered to represent Hilbert's view concerning finitism; see
%Simpson \cite{Si09}.  Note that TPC is a~$\Pi^0_2$-sentence, though it
%is unknown whether the TPC can be settled in WKL$_0$.  The provability
%of TPC in WKL$_0$ would imply the existence of a primitive recursive
%function~$f$ such that~$\phi(x,f(x))$ (is true) for the formula~$\phi$
%in~\eqref{e41}.

With regard to the first sense of the term \emph{infinite totality},
note that it is not the twin prime conjecture that presupposes the
existence of infinite totalities, but rather the assumption that there
must be a determinate answer to whether it is true or false, based on
the belief that the entire universe of all natural numbers exists as a
completed infinity somewhere.  Formalists, including Robinson, refrain
from making such assumptions.

\section{Fallacies and misrepresentations}

We document several fallacies and misrepresentations in Erhardt's
text.

\subsection{A logical fallacy}
\label{s4}

In \emph{Formalism\;64}, Robinson presents his disagreement with his
Platonist opponents in the form of an Alice/Bob-type exchange
{\cite[p.\;231--232]{Ro65}.}
Erhardt misrepresents Robinson's Alice/Bob-type presentation by
attributing Alice's position to Robinson himself, and failing to
mention the more substantive of the two Alice/Bob issues.  To
elaborate, on page 231, Robinson mentions two objections that Alice
(the opponent) might raise:
\begin{enumerate}
  \item[(i)] the `superior brain' argument (which Robinson quickly
    dismisses), and
  \item[(ii)] the `postulation of physical/Platonist infinity'
    argument.
\end{enumerate}    
Alice's `superior brain' argument\;(i) involves the idea that, while
Robinson's brain may have limitations ruling out a ``clear conception
of all sorts of infinite totalities,'' his opponent may possess a
superior conception of the said totalities.

Objection (ii) involves a postulation of infinite totalities in either
the physical or a platonic realm.%
\footnote{For the sake of completeness, we reproduce the full passage
from Robinson: ``An opponent to my position might put forward the
following arguments.

(i) He might say that I am unable to grasp the idea of an actual
infinite totality merely because my brain suffers from a peculiar
limitation.  He might argue that he, on the contrary, has a clear
conception of all sorts of infinite totalities or, at the very least,
of the totality of natural numbers.

(ii) Alternatively, my opponent may concede that he, also, is unable
to grasp the idea of an infinite totality. But he may say that this
does not in any way prove that infinite totalities do not exist. In
order to show that they do exist he may appeal to the physical
world. Or, if he does not wish to, or feels that he cannot, appeal to
the physical world, he may affirm the existence of a platonic world
which contains infinities of all sorts, or of some sort.''
\cite[p.\;231]{Ro65}.}

Robinson quickly dismisses the first objection, pointing out that the
postulation of such superiority ``does not permit any further direct
debate of the issue,'' (Ibid.)
and proceeds to comment on the second objection.  Neither the noun
\emph{psychology} nor the adjective \emph{psychological} occurs in
Robinson's discussion of either objection, or for that matter
elsewhere in his text \emph{Formalism 64}.

Surprisingly, Erhardt attributes Alice's argument\;(i) to Robinson
himself, and misleadingly presents the ``no further direct debate"
comment as Robinson's final word in the issue.  Writes Erhardt:
\begin{quote}
While Robinson has no difficulty \emph{grasping} the possibility of
arbitrarily large finite sets, which motivates his acceptance of the
natural numbers taken individually, he \emph{reports} that it is
simply a \emph{psychological} fact that he cannot grasp infinitary
objects, obviating the possibility of finding common ground with
disputants.%
\footnote{\cite[p.\;431]{Er25}; emphasis on `grasping' in the
original; emphasis on `reports' and `psychological' added.}
\end{quote}
But Robinson `reports' no such thing.  As we already mentioned, the
`psychological' thing is a fabrication.

\subsection{Robinson's letter to G\"odel}

Erhardt's discussion of Robinson's letter to G\"odel contains further
misrepresentations:
\begin{quote}
In addressing a claim that someone with G\"odel’s outlook might
make--that he can fathom the infinite even if Robinson himself
cannot--Robinson states that this objection `does not permit any
further direct debate of the issue'.  (Ibid.)
\end{quote}
Contrary to what Erhardt's quotation (taken out of context) suggests,
Robinson is perfectly willing to engage in a debate - of the pertinent
objection~(ii) (see Section~\ref{s4}).

In his footnote 4, Erhardt similarly misrepresents Robinson's comment
``Hier stehe ich, ich kann nicht anders''
%Luther
in his letter to G\"odel,%
\footnote{This part of Robinson's letter \cite{Ro73} was not
reproduced in \cite{Go03}.}
by quoting it out of context.  Robinson's point is not to concede
purported brain limitations but rather to reject G\"odel's realist
assumptions.  Moreover, Robinson gives his reasons for a reluctance to
accept G\"odel's realist views concerning the hyperreals:
\begin{quote}
[T]he present evidence for the uniqueness of a non standard
$\omega$-ultrapower of the reals is not strong.%
\footnote{\cite{Ro73}.  Recall that hyperreal fields are often
constructed as a quotient~$\mathbb R^{\omega}\!/\mathcal U$
where~$\mathcal U$ is a nonprincipal ultrafilter on~$\omega$.}
\end{quote}

%Is there a name for the logical fallacy of attributing the opponent's
%position (in the context of an Alice/Bob debate) to the writer
%(Robinson) himself?

Erhardt's attribution of Alice's position to Robinson constitutes a
logical fallacy%
\footnote{Robinson is not the only mathematician whose position is
misrepresented by Erhardt.  Erhardt claims that Edward Nelson
``believe[d] there is a largest number'' \cite[p.\;433]{Er25}.  To
support his claim, Erhardt cites four texts by Nelson.  We examined
all four, including Nelson's publication \cite{Ne11}, but found no
evidence of such a belief as claimed by Erhardt.

Erhardt's claim in \cite[p.\;433 note 7]{Er25} that Nelson ``believed
Peano Arithmetic was inconsistent'' is similarly inaccurate.  Rather,
Nelson \emph{suspected} that PA was inconsistent.  When he thought he
had a proof he obviously thought it was inconsistent, but when Tao
found an error in the proof, Nelson immediately accepted~it.}
and does not inspire confidence in his analysis.

\subsection{Verbal excesses from Kreisel to Erhardt}

Similar remarks apply to certain verbal excesses, such as Erhardt's
claim that
\begin{quote}
Robinson has already admitted that such totalities are meaningless,
and without an argument defending the use of these languages,
formalism becomes an \emph{inconsistent} project {\cite[p.\;435;
    emphasis added]{Er25}.}
\end{quote}
Erhardt does not explain how exactly Formalism would become an
`inconsistent' project in the absence of an argument defending the use
of infinite totalities.  Where is the inconsistency?  His claim
appears to be untenable.

Kreisel published his ``Observations on popular discussions of
foundations" in 1971 \cite{Kr71}.  Reviewer James D. Halpern for
Mathematical Reviews described Kreisel's ``Observations'' as ``an
unrestrained attack on P.~J.~Cohen {\ldots}\;and on A. Robinson'' and
noted that
\begin{quote}
Most readers interested in foundations will probably find the
previously mentioned papers of Cohen and Robinson profitable reading,
the disparagement of these in the paper under review notwithstanding.
{\cite{Ha71}}
\end{quote}
As noted by Robinson,
\begin{quote}
Bernays, in an article published in \emph{Dialectica}, criticized my
attitude in his usual gentle manner, while Kreisel had stated his
disagreement with me previously, also in his usual manner.
{\cite[p.\;515]{Ro73b}}
\end{quote}
The reference is respectively to \cite{Be71} and \cite{Kr71}.  In a
critical response to Robinson's text \cite{Ro69}, Bernays claimed that
``there is no fundamental obstacle to attributing objectivity
\emph{sui generis} to mathematical objects.''  {\cite[p.\,178;
    translation ours]{Be71}.}
Robinson appears to have responded in 1975 by including such
`objectivity' of infinitary entities alongside `reference' as claims
that a Formalist would reject; see e.g., \cite[p.\;49]{Ro75}.

In his note 3, Erhardt seeks to distance himself from Kreisel's verbal
excesses, and describes Kreisel's attack as `bias' and `prejudice'.
It is therefore disappointing to find Erhardt himself engaging in
regrettable verbal excesses.

\subsection{Whose \emph{game} was it anyway?}
\label{s5b}

Erhardt lodges the following claim concerning Robinson's view of
mathematics:
\begin{quote}
``He explains non-finitary mathematics as a collection of
  `uninterpretable games with symbols' (Robinson 1969a, p.\;47).''
  \cite[p.\,439]{Er25}
\end{quote}
Looking up the original, one finds that Erhardt has applied a
technique that is already familiar from Section~\ref{s4}: he
attributes Robinson's opponent's position (in this case, `the
intuitionist') to Robinson himself.  For the sake of completeness, we
reproduce the full passage from Robinson's text \emph{From a
formalist's point of view}:
\begin{quote}
[T]he intuitionist may believe that the classical mathematican,
whatever his underlying philosophy, is wasting his time in developing
uninterpretable games with symbols, the sin of the formalist being the
greater because he does so deliberately and consciously.
{\cite[p.\;47]{Ro69}}
\end{quote}
It emerges that, according to Robinson, it is the \emph{intuitionist}
(not Robinson himself) who is wont to make pejorative remarks about
classical mathematicians and formalists allegedly developing
`uninterpretable games with symbols.'%
\footnote{See note~\ref{f11b} for Bishop's comments in this vein.}
Erhardt has again misrepresented Robinson's position.

%\begin{comment}

%{\it

\subsection{Platonism and consistency}

Erhardt opens his text with the following sweeping claim concerning
alleged realist attitudes among mathematicians:
\begin{quote}
[T]he mathematician tacitly adopts a realist outlook, imagining
herself to investigate an \emph{independent realm} of mathematical
objects {\ldots}~about which she can discern objective truths.  This
assumption is crucial to the way mathematicians engage with their
subject.  \cite[pp.\;429--430; emphasis added]{Er25}
\end{quote}
Erhardt's `independent realm' becomes an `independent reality' by page
444 %15. not sure if this was corrected in the second galleys
%
%\footnote{\cite[p.\,15]{Er25}.}
%
Postulating an `independent realm of mathematical objects' is the gist
of a Platonist position.  Thus Erhardt postulates that being a
Platonist is `crucial' for the working mathematician.  Or is it?  Sir
Michael Atiyah, the recipient of both a Fields Medal and an Abel Prize
who was evidently not unsuccessful in ``engaging with his subject",
had the following to say concerning independent realms:
\begin{quote}
The idea that there is a pure world of mathematical objects (and
perhaps other ideal objects) totally divorced from our experience,
which somehow exists by itself is obviously inherent nonsense.
{\cite[p.\;38]{At06}.}
\end{quote}
Erhardt's opening remarks set the tone (and level of seriousness) for
his 20-page text.  In the context of a discussion of independence
results (such as the Continuum Hypothesis), Erhardt claims that
\begin{quote}
[Robinson] \emph{does not consider} the possibility that even if there
are objective answers, they may be beyond our ability to know.
{\cite[p.\;432; emphasis added]{Er25}}
\end{quote}
But Robinson does in fact consider such a possibility, in the
following terms:
\begin{quote}
If~$X$ is the continuum hypothesis, {\ldots}, then we know from the
complementary results of K. G\"odel and P. Cohen that both~$X$ and
non-$X$ are compatible with all known ``natural'' assumptions
regarding the universe of sets (to use platonic language).  While this
suggests to the formalist that the entire notion of the universe of
sets is meaningless (in the sense indicated by our first principle)
\emph{the platonist merely concludes that the basic and commonly
accepted properties of the universe of sets which are known to us at
present are insufficient to decide the continuum hypothesis one way or
the other}.  He will maintain, in this and similar cases, that at any
rate only one of the alternatives that offer themselves is the correct
one, i.e., is in agreement with the truth.  {\cite[p.\;232; emphasis
    added]{Ro65}}
\end{quote}
Robinson then goes on to analyze the difficulties of the Platonist
position with regard to independence results.  Thus, Erhardt's claim
is factually incorrect.

On page 444,
%15 not sure if this was corrected in the second galleys
Erhardt repeats the unfounded accusation aimed at
Robinson, of ``deflating mathematics to a game''
{\cite[p.\;444]{Er25}.}
There follows a curious paragraph on page 15 in Erhardt concerning
Platonism and the problem of consistency:
\begin{quote}
All of this constitutes Robinson giving up on a serious commitment to
his finitism.  This is not to say that he should have pursued a
\emph{definitive argument against platonism}, but that there are more
and less philosophically sophisticated ways of reconciling
mathematical practice with the assertion that we cannot grasp infinite
collections.  (Ibid.; emphasis added)
\end{quote}
Does Erhardt believe that Platonism provides a way out for the problem
of consistency?  The paragraph continues:
\begin{quote}
As it stands, suggesting that mainstream mathematics continue is
tantamount to an admission, despite a reliance on infinitary notions,
that Robinson believes it is consistent.
\end{quote}
This is an allusion to G\"odel's second incompleteness result,
asserting the impossibility of proving Con(PA) (as formalized by
G\"odel) within Peano Arithmetic itself (and similar results for
stronger systems):
\begin{quote}
If not--for instance, if nonstandard analysis led to
contradiction--then anything could be proven from it.
\end{quote}
But such a conclusion would surely hold for just about any piece of
mathematics, and of any philosophical stance.  Or would it?  Erhardt
continues:
\begin{quote}
Robinson does not venture to explain why he takes infinitary
mathematics to be consistent, though a compelling answer might have
provided a rationale for accepting the finitary results it produces.
In the absence of such an explanation, relying upon infinitary
mathematics is itself a source of \emph{great risk}, for there is no
assurance that the derived concrete implications are valid.
{\cite[p.\;444; emphasis added]{Er25}}
\end{quote}
Erhardt holds that relying on infinitary mathematics is a source of
great risk.  The paragraph suggests that Erhardt believes that
Platonism is somehow capable of eluding such a `risk'.  The reasoning
presumably is that if a Platonic realm of mathematics exists, of
necessity it could be neither contradictory nor undecidable, as
Erhardt goes on to comment on issues of faith:
\begin{quote}
Robinson is here relying on \emph{faith} {\ldots}~if not about the
universe of sets, then about basic facts concerning the natural
numbers.  This hunch is not a sufficient warrant to stand in place of
a satisfying mathematical result, like the completion of Hilbert’s
program, and so his instrumentalism is unjustified.  (Ibid.; emphasis
added)
\end{quote}
We leave it to the reader to judge which position involves a greater
leap of faith: that of relying on infinitary mathematics without a
commitment that it possesses a \emph{reference} (in the sense
explained in Section~\ref{s1}), or that of placing one's faith in a
Platonic heaven to ensure the consistency of infinitary mathematics.

Exactly halfway through his text (on page 10 out of 20), Erhardt
finally admits that he is familiar with the philosophical concept of
\emph{not referring}, when he describes Robinson's approach to
non-finitary mathematics as
\begin{quote}
[An] activity that demands no ontological commitments because so many
of the terms in its expressions do not refer.
\\ {\cite[p.\;439]{Er25}.}
\end{quote}
Note that Erhardt's clause ``do not refer'' occurs nowhere in
Robinson's texts, and therefore constitutes Erhardt's own paraphrase
of Robinson's position.  This particular claim of Erhardt's happens to
be an accurate description of Robinson's position concerning the lack
of \emph{reference} (in the sense detailed in Section~\ref{s1}),%
\footnote{To be sure, symbols such as~$\mathbb N$ and~$\mathbb R$ can
\emph{refer} to the corresponding mathematical entities in the context
of a traditional set theory such as ZF.  However, Robinson
specifically speaks of absence of \emph{reference} to entities in the
physical or a putative Platonic realm.}
unlike Erhardt's claims regarding Robinson's purported finitism.

\subsection{Grand justification wherefore?}

Erhardt finds fault with Robinson's strategy as follows:
\begin{quote}
[Robinson] neglects to pursue a \emph{grand justification} for
infinitary mathematics by--in lieu of attempting to offer `a
satisfactory intellectual motivation'--deflating mathematics to a
game.  {\cite[p.\,444; emphasis added]{Er25}.}
\end{quote}
The fact that Erhardt's \emph{game} claim involves a misattribution
has already been addressed: Robinson didn't say it; rather, he
mentioned the intuitionist's use of this pejorative term to criticize
both classical mathematicians and formalists (see Section~\ref{s5b}).
Let us now examine Erhardt's claimed `neglect' of `grand
justification'.  The fact is not merely that Robinson had no such
intention, but that moreover such a project would have been contrary
to his philosophical stance.  Erhardt fails to appreciate the
pragmatic nature of Robinson's philosophical stance, as we now
elaborate.
 
In his 1969 text \cite{Ro69}, Robinson amplifies his list of two items
that he already presented in his \emph{Formalism\;64}: (i)~`infinite
totalities' lack reference (in the sense explained in
Section~\ref{s1}), and
(ii) \emph{useful fictions}: a mathematician
%use specifically singular
should continue acting as if they really existed.

In 1969 Robinson adds a third item: (iii)~\emph{against bottomless
pragmatism}: there is an identifiable core of logical thought that is
common to all of mathematics.  Even more importantly, Robinson uses
these three items to eliminate what he diagnoses as philosophical
dead-ends: items (i), (ii), and (iii) eliminate the positions of,
respectively, the platonist, the intuitionist, and the logical
positivist.

It emerges that Robinson is decidedly not in the business of
developing `grand justification' schemes.  He is merely ruling out
what seem to him to be common philosophical misconceptions, to avoid
ending up in a philosophical dead-end.  To reproach Robinson for not
pursuing a grand justification scheme is to miss entirely the spirit
of Robinson's philosophical posture.

%Quine was not a logical positivist, but like them he subscribed to a
%view that Robinson viewed as nihilist, namely that all, but \emph{all}
%of mathematics is a conventional construct.  From this perspective,
%Erhardt's repeated attempts to prescribe a Quinean framework for
%Robinson%
%
%\footnote{See Erhardt \cite{Er25} note 2; page 10; page 16.}
%
%are inappropriate.

Erhardt claims that
\begin{quote}
While Robinson repudiates actual infinity, the mathematics of his
career revolves around it.  This incongruity mirrors the move from~(i)
to~(ii);
%qua finitist, he must find a way to justify his work, for Gödel’s
%results suggested that finitary consistency proofs of the sort
%envisioned by Hilbert would be insufficient to legitimize the
%continuation of mainstream mathematics for the finitist.
{\ldots}~Robinson \emph{says the same}, noting that \emph{his}
position suffers a serious drawback because `the gap due to the
absence of consistency proofs for the major mathematical theories
appears to be inevitable and we have learned to live with~it'.
\cite[p.\;439; emphasis added]{Er25}
\end{quote}
But there is no `incongruity' here, and Robinson certainly did
\emph{not} ``say the same'' as Erhardt claims.  This is yet another
misrepresentation of Robinson's position, involving a quotation out of
context.  The full quotation makes it clear that Robinson viewed `the
absence of consistency proofs' as a difficulty common to \emph{all}
schools in the Philosophy of Mathematics, and not merely `his
position' as Erhardt claims.  Thus, Robinson wrote:
\begin{quote}
{\ldots}\;\emph{all} known positions in the Philosophy of Mathematics,
including my own, still involve serious gaps and difficulties.  Among
these, the gap due to the absence of consistency proofs for the major
mathematical theories appears to be inevitable and we have learned to
live with it. {\cite[p.\;236; emphasis added]{Ro65}}
\end{quote}
The inevitability of such a gap has been challenged in \cite{Ar25}.
Some logicians are skeptical about the philosophical implications of
Artemov's result.  Regardless of the outcome of that particular
debate, there seems to be little reason to assume that incompleteness
is less of a problem for Platonism than for Formalism, and Robinson
certainly did not claim such a thing.

\section{Types of standardness, realism, knowledge and truth}

We examine some issues related to types of standardness, as well as
some realists' take on the relation of knowledge and truth.  We take a
closer look at some conflations of terms that lead Erhardt to dubious
conclusions about Robinson's work.  We also examine whether the
realism proposed by Erhardt and Gaifman has more convincing answers
than Formalism when it comes to questions of mathematical knowledge
and applicability.

%end \it
%}
%\end{comment}

\subsection{Meanings of Standardness}
\label{s8}

As already mentioned at the end of Section~\ref{s2}, Erhardt appears
to have difficulty keeping apart the generic and the technical meaning
of terms like \emph{meaningless} and \emph{to~grasp}.  Following a
discussion of Skolem's nonstandard models and Robinson's nonstandard
analysis, involving the \emph{technical} distinction between standard
and nonstandard numbers, Erhardt writes:
\begin{quote}
``Adding to the irony--and \emph{irresolvable tension} in Robinson’s
  work and philosophy--is that, due to its predication on the
  transfinite, Robinson himself says that `the entire notion of
  \emph{standardness} must be meaningless [\ldots]' (Robinson
  1964/1979b, p. 242).'' \cite[p.\;442 note 33; emphasis added]{Er25}
\end{quote}
The irony, according to Erhardt, is that Robinson himself talks
`standard and nonstandard', and then comes out swinging against `the
entire notion of standardness'.  However, the context - not clarified
by Erhardt - of Robinson's comment on `the entire notion of
standardness' is a discussion of the Platonist idea of a standard or
intended model of arithmetic or set theory, as analyzed in
Section~\ref{s13}.  To reject such an idea, Robinson used the term
`standard' in its generic sense, whereas Erhardt misleadingly presents
it as if what is involved is the technical sense used in nonstandard
analysis.

In more detail, theories of nonstandard analysis exploit a
mathematical concept called the standardness predicate, which
distinguishes between standard and nonstandard numbers (and more
general sets).  Since this distinction can be thought of as a
formalisation of Leibniz's distinction between assignable and
inassignable numbers,%
\footnote{See e.g., the publications \cite{13f}; \cite{23f};
\cite{24b}; \cite{24d}, \cite{25d}.}
the predicate could be referred to also as the \emph{assignability
predicate}.  Meanwhile, Platonists believe in (and Robinson rejects)
the existence of standard models a.k.a. \emph{intended
interpretations} of $\N$, $\R$, and other mathematical entities.  When
the issue is translated into the terminology of \emph{the
assignability predicate} and \emph{the intended interpretation}, it is
obvious that there is no connection between them, contrary to the
impression Erhardt seeks to create.

%\begin{comment}
%{\it
  
\subsection{Resolving tensions}

It is worth examining Erhardt's claim concerning an alleged
`irresolvable tension in Robinson's work and philosophy' in more
detail.  Here one needs to distinguish between \emph{doing}
mathematics and \emph{talking about} mathematics.

(1) The question whether~$10000000000000001$ is a prime
%if nis not a power of 2 then 10^n+1 is obviously not prime.  Here
%n=16.
concerns a concrete natural number.  One does not need Peano
Arithmetic, let alone the infinite set~$\N$, for it to be meaningful.

(2) On the other hand, to make sense of the question whether there are
arbitrarily large numbers~$n$ for which~$10^n + 1$ is prime does
require some theory of arithmetic, and perhaps (if Erhardt's reading
of Quine is right) also an infinite set~$\N$.

(3) The notions of formal logic such as \emph{variable}, \emph{term},
\emph{formula} and \emph{proof} are defined by recursive rules similar
to the definition of natural numbers.  We therefore have an analogous
distinction to that between items (1) and (2) above, as follows.

(4) When we write down (more/less formal) actual proofs in some area
of mathematics, these proofs are concrete objects analogous to
concrete numbers such as~$10000000000000001$.  We do not need the
collection of all variables, terms, formulas or proofs to verify that
a given concrete proof really is a correct proof.  This is what
formalists like Robinson are doing when proving results in some area
of \emph{mathematics}, including aeronautics, they are working in; no
actual infinity is necessarily involved here.

(5) Formalists can also be interested e.g., in determining \emph{what
is provable}.  Then they are doing metamathematics, i.e., the area
they are working in is \emph{logic}.  For this purpose one will likely
want to exploit the kind of infinite collections listed in item (3).

Thus, contrary to Erhardt's claim, there is no tension between
Robinson's work in aeronautics (mathematics), on the one hand, and his
work in - or philosophy of - metamathematics (logic), on the other.
One could make the following additional points.

(6) When mathematicians practice some area of mathematics, they prove
theorems.  Such theorems (for example, Wiles’s proof of Fermat's Last
Theorem) are concrete objects.  No theory of proofs, and in particular
no commitment to infinite collections, is needed to generate such
proofs.

(7) A formalist mathematician may feel compelled (but, arguably, is
not required \emph{qua} philosopher) to study mathematical practice
mathematically.  If and when doing that, they practice mathematics in
another field: logic, proof theory, model theory{\ldots} In such
practice, they may find (referenceless) infinite collections of
variables, formulas or proofs useful.  The remark of item (6) applies
to such practice, as well.  There is no need for actual infinity in
order to give a proof of and find useful G\"odel’s Incompleteness
theorem.  Thus even in logic there are results that do not depend on a
commitment to any actual infinity.

In sum, there is neither contradiction nor tension between a
formalist’s practice and theory (philosophy).

\subsection{Varieties of realism}
\label{s5}

Not all mathematical realists hold identical views.  We will examine
Erhardt's position in a framework proposed by Hamkins.  In an
influential 2012 article, Hamkins provides the following useful
perspective:
\begin{quote}
The \emph{multiverse} view in set theory
%  introduced and argued for in this article,
{\ldots}~is the view that there are many distinct concepts of set,
each instantiated in a corresponding set-theoretic universe.  The
\emph{universe} view, in contrast, asserts that there is an absolute
background set concept, with a corresponding absolute set-theoretic
universe in which every set-theoretic question has a definite answer.
The multiverse position, I argue, explains our experience with the
enormous range of set-theoretic possibilities, a phenomenon that
challenges the universe view.  {\cite[Abstract]{Ha12}; see also
  \cite[Section~8.14]{Ha21}}
\end{quote}
Hamkins goes on to state his position with regard to the Continuum
Hypothesis (CH):
\begin{quote}
In particular, I argue that the continuum hypothesis is settled on the
multiverse view by our extensive knowledge about how it behaves in the
multiverse, and as a result it can no longer be settled in the manner
formerly hoped for.  (Ibid.)
\end{quote}
In more detail, Hamkins views CH as a ``switch'' that can be turned
on and off as it were, by passing to a suitable forcing extension of
the current instance of a set-theoretic universe.

Hamkins considers himself a realist of the multiverse, thereby
discarding the assumption of the \emph{uniqueness} of a set-theoretic
universe as sometimes held by other realists.

We will now examine Erhardt's stance, in the context of Hamkins'
dichotomy.  Erhardt writes:
\begin{quote}
[T]he mathematician tacitly adopts a realist outlook, imagining
herself to investigate an independent realm of mathematical
objects--one that includes infinite totalities like the natural
numbers--about which she can discern objective truths.
\\ {\cite[p.\;430]{Er25}.}
\end{quote}
While this passage suggests that Erhardt's sympathies are with the
realist, it does not indicate whether or not this constitutes realism
about a \emph{unique} set-theoretic universe.  Reading further, we
find the following comment concerning Zermelo--Fraenkel set theory
with the Axiom of Choice (ZFC) and the CH:
\begin{quote}
Consider Cantor’s continuum problem: Is there an intermediate
cardinality between the naturals and the reals?  If you are a realist
about infinite totalities, then the continuum problem is an obvious
question to ask. However, ZFC+CH and ZFC+CH~[sic]%
\footnote{The second occurrence of ZFC+CH is a typo and should be
``ZFC+$\neg$\,CH''.}
are both consistent, which Robinson takes as implying that there must
be no universe of sets. A platonist will respond that we have simply
not yet found the axioms stating the basic properties of the universe
of sets sufficient to decide CH one way or the other.  \emph{This is a
valid objection}, but it does not sway Robinson.  {\cite[p.\;432 note
    5; emphasis added]{Er25}}
\end{quote}
Here Erhardt is sympathetic to the idea that CH ought to possess a
definite truth value yet to be determined, and considers such
reasoning as `valid'.  Erhardt's position concerning CH is clearly at
odds with Hamkins' multiverse-realism, and specifically with Hamkins'
position with regard to CH (namely, that it need \emph{not} possess a
definite truth value yet to be determined).  Our tentative conclusion
is that Erhardt is a universe-realist.%
\footnote{Kreisel asks: ``Can we extend the language of set theory to
{\ldots}\;decide CH by means of axioms (which are evident for the
intended interpretation) of the extended language?''
\cite[pp.\;196--197]{Kr71}.  Thus far, evidence in favor of an
affirmative answer has been weak.}
As argued above, this may color his evaluation of Robinson's
philosophical position.

\subsection{Applicability of mathematics}

The far-reaching applicability of mathematics in diverse areas is
certainly a famous problem without easy answers.  However, when
Erhardt claims that
\begin{quote}
[Robinson's] foundational position
%he inherited from David Hilbert
undermines not only the use of nonstandard analysis, but also
Robinson's considerable corpus of pre-logic contributions%
\footnote{See note \ref{f11}.}
to the field in such diverse areas as differential equations and
aeronautics {\cite[Abstract]{Er25}}
\end{quote}
one may well wonder how exactly the postulation of a Platonic
set-theoretic universe would help explain the applicability of
mathematics in ``such diverse areas as differential equations and
aeronautics.''  In other words, how would the putative existence of a
mind-independent realm of mathematical abstracta help explain the
ability of an embodied brain to apply the said abstracta to the
solution of concrete problems of airplane wing design
(cf.~\cite{Ro56})?  It would seem that the difficulty involved is no
less challenging than the applicability of what Leibniz already
referred to as \emph{useful fictions}.%
\footnote{See \cite{13f}; \cite{24b}; \cite{24d}; \cite{25d}.}

Robinson proposed an apt diagnosis of the Platonist's aversion to the
Formalist's position:
\begin{quote}
[I]t is perhaps natural that a mathematician should resent suggestions
which deprive him of the comforting feeling that he, like the
physicist or biologist, spends his life in the exploration of some
form of reality.  {\cite[p.\;46]{Ro69}.}
\end{quote}
As far as the criterion of the applicability of mathematics in the
natural sciences, Robinson explains:
\begin{quote}
[T]his criterion does not imply that the terms of the theory should be
interpretable directly and in detail.  It is sufficient that we should
have rules which tell how to apply certain relevant parts of our
theory to the empirical world.  {\cite[p.\;234]{Ro65}.}
\end{quote}
Thus, infinitary terms may be \emph{referenceless}, but their
applications can nonetheless be \emph{meaningful}.

%end \it
%}
%\end{comment}

\subsection{Gaifman's realism}
\label{s12}

\cite{Ga12} claims to use the term \emph{realism} in the sense of
%Shapiros
\emph{realism in truth-value}.  In relation to the natural numbers,
this means assuming that every PA sentence is a factual statement,
i.e., has an objective truth value (even if it is undecidable in PA
itself).  As Robinson pointed out almost half a century earlier,
\begin{quote}
[A]s a matter of empirical fact the platonists believe in the
objective truth of mathematical theorems \emph{because} they believe
in the objective existence of mathematical entities.  {\cite[p.\;230;
    emphasis in the original]{Ro65}}
\end{quote}
There is internal evidence in Gaifman's article that such is indeed
his position, as illustrated by his snide comment about Hamkins and
Robinson:
\begin{quote}
Another possibility has emerged from the multiverse conception
proposed by Hamkins (2011).  On this view, the model of natural
numbers depends on the set-theoretic universe containing the model.
What it comes to is that we have a clear enough conception of models
of PA, or of some extensions of PA, but we have no clear distinction
between standard and non-standard models.  Perhaps Abraham Robinson,
the founder of non-standard mathematics, who was not a set
theoretician, would agree to that.
{\cite[p.\;489]{Ga12}}
\end{quote}
Gaifman's professed truth-value realism appears to be consistent with
the belief that the `standard or intended model of arithmetic' is a
meaningful notion.  Recall that Robinson rejects such a notion (see
Section~\ref{s13}).

Note that ``$\phi$ is provable'' and ``$T$ is consistent'' are
respectively~$\Sigma^ 0_1$ and~$\Pi^0_1$ sentences (when encoded in
PA).  Gaifman apparently holds that one commits oneself to realism in
truth value (in relation to the natural numbers) when one uses the
terms \emph{provable} and \emph{consistent} in metamathematics.  He
seems to find it philosophically unacceptable to assume that
metamathematical statements about the provability of a sentence or the
consistency of a theory could be undecidable.  Meanwhile, for a
formalist practicing finitistic metamathematics in Primitive Recursive
Arithmetic (PRA) such a posture toward undecidability is a natural
consequence.  The “finitistic semantics” of PRA is naturally
incomplete; {see \cite{Ta81}.}
A~$\Pi^0_1$ sentence is `true' if it can be proved by induction; it is
`false' if one can provide a counterexample.  A $\Sigma^0_1$ sentence
is `true' if one can provide an example; it is `false' if one can
provide an inductive proof of the corresponding $\Pi^0_1$ sentence
expressing its negation; see further in \cite{Si09}.

Similarly to intuitionism, \emph{truth} in this semantics is defined
by provability, and \emph{falsity} by refutability.  Consequently, in
metamathematics, on the basis of potential infinity, there is no
\emph{tertium non datur}. Such a position is unacceptable for Gaifman.

Such a metamathematics can nevertheless prove meaningful statements.
For instance, as explained in \cite{21e} the conservativity of the
theory SPOT over ZF can be proved in WKL$_0$ (which is not
finitistic).  But this conservativity result is a~$\Pi^0_2$ sentence,
and according to the results of reverse mathematics \cite{Si09}, there
is therefore a proof in PRA, where the result is formulated
as~$\phi(m, f(m))$ with a primitive recursive function.  Such a
formulation is finitistic.

The fact that realist positions (like Gaifman's) make it harder to
appreciate the hidden resources of the theory of the real number line
was noted by Massas in the following terms:
\begin{quote}
[T]he conservativity results in [Hrbacek and Katz]%
\footnote{Massas is referring to \cite{21e}.}
would likely fall short of convincing anyone who thinks that any
existence claim regarding Robinsonian infinitesimals is simply false.
\\ \cite[p.\;263]{Ma24}
\end{quote}

Paul Cohen seemed to regard it as a prerequisite for his progress in
foundational research to abandon the idea that there is a unique
(intended) interpretaion of $\mathcal P(\N)$ that conforms to our
intuition:
\begin{quote}
I can assure that, in my own work, one of the most difficult parts of
proving independence results was to overcome the psychological fear of
thinking about the existence of various models of set theory as being
natural objects in mathematics about which one could use natural
mathematical intuition.  {\cite[p.\,1072]{Co02}}
\end{quote}

\subsection{Knowledge-transcendent truth?}

With regard to the relation between independence results and
mathematical realism, Gaifman mentions `knowledge-transcendent truth':
\begin{quote}
If an independence result indicates that some mathematical truths
outstrip our capacities of knowing them, then it points to
knowledge-transcendent truth, hence to realism.
{\cite[p.\;498]{Ga12}}
\end{quote}
Objectively speaking, the implication would seem to go the other way
around: it is only \emph{if} one wishes to defend realism that one may
be led to postulate knowledge-transcendent truth.  Such
`knowledge-transcendent truth' is mentioned twice in the article (the
second occurrence is in the concluding section).  It emerges that when
Gaifman claims that something points toward realism, what is he has in
mind is that it points toward \emph{transcendent truth}.  However,
Gaifman's argument appears to be circular.  Namely, the starting
assumption that ``independence results indicate that some truths
outstrip our capacity'' itself depends on a realist standpoint.  An
anti-realist has little reason to assume that independence results
would indicate such a thing.  So the full Gaifmanian circularity would
run as follows:
\[
\begin{aligned}
\text{\small realism} &\Rightarrow \text{\small independence results
  indicate that truths outpace our capacities} \\&\Rightarrow
\text{\small there exist transcendent truths} \\&\Rightarrow
\text{\small realism}.
\end{aligned}
\]
Gaifman claims further:
\begin{quote}
I am mostly concerned with independence results for first-order
arithmetical statements.  My general view is that these results do not
affect our conception of the standard natural numbers.  Hence, they
reinforce any realism that is based on this conception.
{\cite[p.\;509]{Ga12}}
\end{quote}
G\"odel's second incompleteness theorem as applied to PA is a classical
arithmetic independence result.  Gaifman appears to claim that such a
result would ``reinforce realism.''  Accordingly, the independence of
CON(PA) of PA would reinforce realism.  One wonders what Gaifman would
make of Artemov's recent proof of CON$^S$(PA) within PA.%
\footnote{\label{Art}\cite{Ar25} proposed an alternative formalisation
Con$^S$(PA) of consistency following Hilbert.  Giaquinto quotes
Robinson in \cite[p.\;239]{Ro65} to the effect that ``I cannot see at
this time how a form of reasoning which attempts to escape the
consequences of G\"odel's second theorem (such as Gentzen's
consistency proof or any other consistency proof for Arithmetic) can
remain strictly finitistic and hence interpreted'' and claims that
``the formalist foundational programme is irredeemable''
\cite[p.\,129]{Gi83}.  Artemov's result provides a possible way of
escaping the consequences of G\"odel's second theorem, and therefore
possibly undermines Giaquinto's rash conclusion.  Similar remarks
apply to Giaquinto's surprising claim that ``the resuscitation of
formalism by Robinson and Cohen (following Cohen's proof of the
independence of the Continuum Hypothesis) is untenable and barren''
\cite[p.\,120]{Gi83}.}
If even independence reinforces realism, lack of independence (in the
sense of Artemov's result) would seem to also reinforce realism.  This
would suggest that Gaifman's variety of realism is unfalsifiable.

All known undecidable first-order arithmetic sentences have a
preferred truth value.  For example, CON(PA) is generally considered
to be true (and is of course provable in a stronger system such as
ZFC), unlike higher-order statements such as CH.  Some mathematicians
view this as an argument in favor of realism about PA.  However, such
an argument does not appear in \cite{Ga12}.

\section{Conclusion}

While Erhardt goes on eventually to discuss Robinson's mathematical
Formalism, his opening presentation of Robinson's position as finitism
is off the mark, and his attempt to present Robinson's opponent's
position as Robinson's own does not inspire sufficient confidence to
give credence to his speculations concerning Robinson's analysis of
potential infinity.

Erhardt's endorsement of a naive realism about a standard model of
ZFC, including his faith in a determinate truth value of CH (yet to be
`discovered'), fails to take into account the recent literature in the
philosophy of the foundations of mathematics, such as the multiverse
perspective of \cite{Ha12}; see also \cite[Section~8.14]{Ha21}.  Due
to his misreading of Robinson's terms \emph{meaningless}, \emph{to
grasp}, and \emph{standard}, Erhardt sees an irony and an unresolvable
tension between Robinson's use of nonstandard models of the real
numbers and his statement that ``the entire notion of standardness
must be meaningless.''  However, there is neither irony nor tension
here, since Robinson is referring to his previous statement that
infinite sets, including both standard and non-standard models, lack
reference.  He did not mean that they are pointless or devoid of
significance.

In his 1969 text, Robinson envisions
\begin{quote}
{\ldots}~a solution which involves the equal acceptance of
\emph{several kinds of set theory}.  At this point, it is natural to
recall the historical development of Euclidean geometry whose analogy
with the recent development of axiomatic set theory is perhaps closer
than is generally appreciated.  {\cite[pp.\;46-47; emphasis
    added]{Ro69}}
\end{quote}
The comparison between the existence of distinct set theories and the
existence of distinct geometries was pursued in more detail by
\cite{Co67}.  Different set theories of this type were indeed
developed shortly after Robinson's passing.  \cite{Hr78} and
\cite{Ne77} independently developed axiomatic approaches to
nonstandard analysis as an alternative to model-theoretic approaches.
In fact, Robinson's visionary insight is arguably more likely to
emerge from a Formalist than a Platonist philosophical stance.%
\footnote{Our conclusion here is at variance with Kreisel's.  Kreisel
claimed the following: ``having accepted the formalist scheme, one is
not merely \emph{unwilling} to appreciate (or study!) further work in
foundations, but simply \emph{unable} to do so without reexamining
one's ideas'' \cite[p.\,191]{Kr71} (emphasis in the original).
However, Hrbacek and Nelson were indeed \emph{able} to appreciate,
study, and make progress in foundations without adopting Platonist
ideas.}

While traditional set theory is formulated in the~$\in$-language, the
axiomatic approaches enrich the language of set theory by means of the
introduction of a one-place predicate `standard' or `st', and are thus
formulated in the st-$\in$-language.  Such theories are conservative
over ZFC.\, Suitable axioms govern the interaction of the new
predicate with the traditional ZFC axioms.  In the axiomatic
approaches, infinitesimals are found within~$\R$ itself (rather than
in an extension, as in the model-theoretic approaches to nonstandard
analysis), and nonstandard integers in~$\N$ itself.  This challenges
Platonist notions about~$\N$ and~$\R$ as determinate (or even
mind-independent) entities.  Robinson's Formalism remains a viable
alternative to mathematical Platonism.

%\section{Competing interests}There are no competing interests.

\section*{Acknowledgments}We are grateful to Karel Hrbacek for insightful comments.

\end{document}